\documentclass{article}
\usepackage{amssymb}
\usepackage{amsbsy}
\usepackage{epsfig}
\newtheorem{theorem}{Theorem}[section]
\newtheorem{corollary}{Corollary}[section]
\newtheorem{lemma}{Lemma}[section]
\newtheorem{example}{Example}[section]
\newtheorem{remark}{Remark}[section]
\newtheorem{proposition}{Proposition}[section]
\newtheorem{definition}{Definition}[section]
\def\pa{\partial}
\def\be{\begin{equation}}
\def\ee{\end{equation}}
\def\bes{$$}
\def\S{\Sigma}
\def\l{\lambda}

\def\q{{\mathfrak q}}
\def\Exp{{\boldsymbol \exp}}
\def\ees{$$}
\def\le{\left}
\def\ri{\right}
\def\bea{\begin{eqnarray}}
\def\eea{\end{eqnarray}}
\def\beaq{\begin{eqnarray}}
\def\eeaq{\end{eqnarray}}
\def\beas{\begin{eqnarray*}}
\def\eeas{\end{eqnarray*}}
\def\k{{\varkappa}}
\def\h{{\mathfrak h}}
\def\D{\nabla}
\def\bl{\begin{lemma}}
\def\el{\end{lemma}}
\def\bp{\begin{proposition}}
\def\ep{\end{proposition}}
\def\bc{\begin{corollary}}
\def\ec{\end{corollary}}
\def\bt{\begin{theorem}}
\def\et{\end{theorem}}
\def\br{\begin{remark}}
\def\er{\end{remark}}

\def\bx{\begin{example}}
\def\ex{\end{example}}

\def\bd{\begin{definition}}
\def\ed{\end{definition}}

\def\be{\begin{equation}}
\def\ee{\end{equation}}

\def\R{{\mathbb R}}
\def\o{\omega}
\def\O{\Omega}
\def\p{\pi}
\def\s{\sigma}
\def\Ho{H^{\o}}
\begin{document}
\begin{flushright}
SISSA 72/99/GEO
\end{flushright}
\vspace{0.2cm}
\begin{center}
\begin{Large}
\textbf{Rigidification of  pseudo--Riemannian
 Manifolds by an Elliptic Equation}\footnote{
Work supported in part by the Natural Sciences and Engineering Research Council
of Canada (NSERC) and the Fonds FCAR du Qu\'ebec.}
\end{Large}\\
\vspace{1.0cm}
\begin{large} {Marco
Bertola}$^{\dagger\ddagger}$\footnote{bertola@crm.umontreal.ca}
 and Daniele Gouthier$^{\star}$\footnote{gouthier@sissa.it}
\end{large}
\\
\bigskip
\begin{small}
$^{\dagger}$ {\em Centre de recherches math\'ematiques,
Universit\'e de Montr\'eal\\ C.~P.~6128, succ. centre ville, Montr\'eal,
Qu\'ebec, Canada H3C 3J7} \\
\smallskip
$^{\ddagger}$ {\em Department of Mathematics and
Statistics, Concordia University\\ 7141 Sherbrooke W., Montr\'eal, Qu\'ebec,
Canada H4B 1R6}\\
$^{\star}${\em SISSA, V. Beirut 2--4, 34014 Trieste, ITALY }\\
\end{small}
\bigskip
\bigskip
%%%%%%%%%%%%%%%%%%%%%%%%%%%%%%%%%%  Abstract %%%%%%%%%%%%%%%%%%%%%%%%%%%%%%%
{\bf Abstract}
\end{center}

 We study the solvability of the equation for the smooth function $\o$, 
 $\Ho=-\k\o\,g$ on a geodesically complete
pseudo--Rie\-man\-nian  manifold $(M,g)$, $\Ho$ being the covariant
 Hessian of $\o$. A similar equation was considered by Obata and Gallot in the
 Riemannian case for positive values of the constant $\varkappa$; the
 result was that the manifold must be the canonical sphere.
 In this generalized setting 
we obtain a range of possibilities,
 depending on the sign of $\k$, the signature
 of the metric and the value of a certain first integral of the equation:
the manifold is shown to be of constant sectional curvature or a
 warped product with suitable factors depending on the cases. 

%%%%%%%%%%%%%%%%%%%%% Paper begins %%%%%%%%%%%%%%%%%%%%%%%%%%%%%%%%%%%%%%%%

%
%\baselineskip 20pt plus 2pt minus 2pt
\section{Introduction.}

Given a Riemannian manifold $(M, g)$, we denote by $\Ho$ the Hessian of the
 smooth real function $\o$.
It is a classical result that the
 existence of a solution of $\Ho=-\o\, g$ constrains the curvature to be equal
 to $1$ \cite{kn:ob}. Such theorem can be generalized  by the equation  
\be
\Ho = -\k\o\,g\ ,\label{magn}
\ee
where now $\k\in \R$ and $(M,g)$ is a pseudo-Riemannian manifold.\\
This equation was considered in the literature 
 only in the Riemannian case and for positive
 values of the constant $\k >0$  by Obata \cite{kn:ob} and Gallot \cite{kn:ga}.
In particular Obata states the theorem (Obata's Theorem): if $(M,g)$ is geodesically complete
 and $\o$ is a nontrivial solution to eq. (\ref{magn}) with $\k>0$, then
 $(M,g)$ is isometric to the canonical sphere of radius $R=\frac
 1{\sqrt{\k}}$.\\ 
Unfortunately the proof of this nice theorem cannot be generalized to
 this context.
The purpose of this paper is twofold; first
 we provide a proof  that works independently of the sign of
 $\k$ in the Riemannian case (the sphere being replaced by a
 negatively curved hyperboloid for $\varkappa<0$).
Secondly, the technique of the proof can be applied to the
 pseudo--Rie\-man\-nian cases as well.\\
For  historical reasons we name the equation (\ref{magn}) ``Obata's
 Equation''  and a solution of it a ``Obata's function''.\par
A related equation was considered by Kerbrat in \cite{kn:ke}; in his case the
 equation is of third order in $\o$, namely 
\beaq
\nabla_Y\Ho(X,Z) = -\k\le(
 2{\rm d}\o(Y)\langle X,Z\rangle+{\rm d}\o(X)\langle Y,Z \rangle +
 \langle X,Y\rangle{\rm d}\o(Z)\ri)\ ,
\eeaq
which is presented and studied in 
the case $\varkappa = -1$; the author then proves a certain rigidification assuming the
 existence of at least two linearly independent solutions or a full
 rigidification if the critical set of the solution is nonempty, a
 result that
 corresponds to our Thm. (\ref{reduction}) in the present context.\\
The relation between these equations is elucidated by Gallot in
\cite{kn:ga}.
 The author  shows that (for $\k>0$) the two equations are
the first two of a sequence of equations $E_n$ of $n$-th order
 equations characteristic of the canonical sphere and he shows  that their
solutions are  the harmonic homogeneous polynomials of degree $n-1$ in
the Euclidean space in which the sphere is embedded: the solutions of
 equation $E_n$ are proper eigenfunctions of the Laplacian.\par
The study of solvability of Obata's equation
depends also on the value of a first integral, namely
\be
\|\O\|^2 + \k\o^2 = \h\ ,\qquad \O:= {\rm grad}(\o) = {\rm
 d}\o^\sharp \label{firstint2}
\ee
As outlined above, 
we generalize Obata's theorem in two directions: first, we allow any
 real value of $\k$ while keeping the assumption on the Riemannian
 signature of the metric.
 The result is that the  manifold is of constant sectional 
 curvature only if $\k\h>0$ ($\h$ being the first integral in
 Eq. \ref{firstint2}), while in all other cases we have only the
 splitting of the metric in a ``warped--product'' metric
(cf. \cite{kn:on}), without any further constraint on the Riemannian 
 curvature.\par
Secondly, we consider the same equation on a pseudo--Rie\-man\-nian
 manifold.
Once more we obtain that the sectional curvature must be a constant in
 the case $\k\h>0$ while for $\k\h<0$ we obtain  only the splitting in
 a warped product metric.
 However we also obtain an ``intermediate''
 rigidification in the case $\k\neq 0$, $\h=0$,
 with a  necessary  condition of asymptotic flatness on the level
 surfaces of $\o$. 
This has no analog in the Riemannian case, where for $\h=0$
 (and {\em a fortiori} $\k<0$) we only have the warped--product
 structure.
In order to construct a nontrivial example, 
in Prop. \ref{counterex} we define a geodesically 
complete manifold supporting such an Obata's 
function by gluing two suitable incomplete manifolds along a lightlike smooth
hypersurface. 
This manifold does  {\em not} have constant curvature in general,
 showing that this is not a necessary condition, and is not globally a
 warped-product.
Another case which is not present in the Riemannian setting is for
 $\h=\k=0$, because in Riemannian signature this would immediately
 imply that $\o\equiv 0$ (using eq.(\ref{firstint2}) and
 (\ref{magn}));
on the contrary, in pseudo--Riemannian curvature the equation has
 nontrivial solutions corresponding to a null Killing vector field
 (Thm. \ref{mainpseudo4}, part (iii)).\par
The paper is organized as follows: in Section (\ref{prelim}) we set some
notation recalling notion and properties of warped products 
and we establish  some preliminary results on geodesic completeness of
pseudo--Rie\-man\-nian warped products.
In   Section (\ref{main}) we state and prove
 the theorem of rigidification for Riemannian manifolds.
Many statements do not rely on the signature and hence can be rephrased
 without any change in the pseudo--Rie\-man\-nian case.
In Section (\ref{mainps}) we extend this theorem to the case of arbitrary
 signature: since the result and the proof 
depend strongly on the relative signs of $\k\,\h$, we
 split the theorem in Thms. \ref{mainpseudo1}, \ref{mainpseudo2},
 \ref{mainpseudo3}, \ref{mainpseudo4} in order to avoid too many case distinctions.
Finally Section (\ref{maximal}) is devoted to the study of a maximal set of
 Obata's functions, where it is shown that if there exist more than one
 solution, then the distribution spanned by their gradients is involutive and 
 foliates the  manifold in submanifold of constant sectional curvature.
\section{Preliminaries}
\label{prelim}
We first recall the notion of  warped product.
Let  $(B, g_B)$ and $(F, g_F)$ be two pseudo--Rie\-man\-nian manifolds. 
Consider the smooth
manifold $M := B \times F$ with the canonical projections denoted by
$\p:M \rightarrow B$ and $\s:M \rightarrow F$. Given an arbitrary
smooth map $\alpha: B \rightarrow \R^+$ we define a 
 (pseudo)--Rie\-man\-nian metric $g=g_{\alpha}$ on $M$
(called {\it warped metric} \cite{kn:on})
\be
g_{\alpha} := \p^* g_B + (\alpha \circ \p)^2 \s^* g_F .
\ee
  Let  $X,Y$ be sections of $\Gamma(\pi^*TB)$ and
$U,V$ of $\Gamma(\sigma^*TF)$ and let $A$ denote the gradient
of $\alpha$.
 The sectional curvature is given in terms of the sectional curvatures of $B$
and $F$ as follows \cite{kn:bg}:
\beaq
K_{XY} = K^B_{XY} ;\qquad
K_{XV} = - {{H^\alpha(X,X)}\over{\alpha ||X||^2}} ;\qquad
K_{UV} = {{K^F_{UV} - ||A ||^2}\over{\alpha^2}} ,
\label{sectional}
\eeaq
where in these formulas the norms are the pseudo--lengths.\par
It is a well-known result that in the Riemannian case 
 a warped product $B\times_\alpha F$ is geodesically complete
iff both factors are and $\alpha>0$ \cite{kn:bo}. 
The pseudo--Rie\-man\-nian case is much less
studied in general.  If the base $B$ is one dimensional we can
establish the following lemmas.
\bl
\label{completeness}
Let $M=B\times_\alpha F$ be a pseudo--Rie\-man\-nian warped product with $B$
(anti)--Rie\-man\-nian and
both factors geodesically complete.\\
If  $\epsilon:= \mathop{\hbox{\rm inf}}_B\alpha>0$ then  $M$
is geodesically complete.
\el
{\bf Proof}. Consider the equation for a geodesic $\gamma(s)$ of type space,
time or light and set correspondingly $C=+1,-1,0$:
 decomposing  the vector $\dot \gamma = X+V=\pi_*(\dot \gamma)+
\sigma_* (\dot \gamma)$  we get (cf. \cite{kn:bo})
\be
\nabla^B_X X= \frac {\|V\|^2}\alpha A\ ,\qquad 
\nabla^F_V V= -2\frac {X(\alpha)}\alpha V\ .
\ee
Computing the rate of change of the norms we get
\be
-\frac d{ds} \le(\|V\|^2\ri)=
\frac {d}{ds} \le(\|X\|^2\ri) = 2\frac {X(\alpha)}\alpha \le(C-\|X\|^2\ri)
\ .
\ee
Integrating the equation once we  obtain 
\be
\|X\|^2 = (C-\|V\|^2) = \frac{{\alpha_0}^2 (\|X_0\|^2-C)}{\alpha^2\circ
\pi(\gamma(s))} +C \leq \frac {{\alpha_0}^2 |\|X_0\|^2-C|}{\epsilon^2} +|C|\ .
\ee
This shows  that the square norm of $X$ is bounded, hence the curve $\pi\circ
\gamma:I\to B$ has finite length for any value of the parameter.
This and geodesic completeness (which, for (anti)--Rie\-man\-nian manifolds is
equivalent to completeness)  proves that the curve $\pi\circ\gamma$ is defined
for any $s\in \R$.
Moreover   the norm of $V$
\be
\|V\|_F^2 = \frac 1 {\alpha^2} \|V\|^2 = \frac
 {C-\|X\|^2}{\alpha^2}
\ee
 is bounded as well since $\|X\|^2$ is.
 Now $V$ is parallel translated because the projection of
$\gamma$ on the second factor is a pregeodesic whose parameter is bounded
by a suitable affine parameter, hence $\sigma\circ \gamma$ is defined for any
 $s\in\R$ as well.\\
This proves that  $M$ is
geodesically complete. $\rule{5pt}{8pt}$\par\vskip 3pt
\bl
\label{incompleteness}
Let $M=\R\times_\alpha \Sigma$ be a pseudo--Rie\-man\-nian manifold of type
$(r,p)$ with metric ${\rm d}s^2 = -{\rm d}t^2 +\alpha^2(t) g_\S$.
 Suppose $\alpha$ is monotonic out of  a bounded set and $\displaystyle{\mathop
{\hbox{\rm inf }}_{t\in \R}\,\alpha(t)=0}$.\\
If $\alpha$ is integrable at either $ +\infty$ or $-\infty$ then $M$ is not
geodesically complete.\\
An inextensible  incomplete geodesic reaches the  spacelike boundary of
$\Sigma$ for a finite value of the affine parameter.
If the base $\R$ is (anti)Riemannian, then the geodesic
 reaches the (spacelike) timelike boundary
of $\Sigma$.
\el
{\bf Proof}.  The proof is based on the previous formulas. Suppose
that the factor $\R$ is of Riemannian signature; if $s$ is the affine
parameter of the geodesic, the equation for the geodesic coordinate
$t$ of $\R$ can be recast into
\be
\frac {\alpha(t) {\rm d}t}{\sqrt{ \alpha^2 (t_0)({\dot t_0}^2+ C) -
C\alpha^2(t)}}  = {\rm d}s\ ,\label{affineparameter}
\ee
where $C=+1,-1,0$ according to the type space, time or light of the
geodesic.\\
For timelike geodesics ($C=-1$) for which $|\dot t_0|>1$ 
we can have unbounded trajectories in the
direction(s) where $\alpha$ goes to zero (namely one or both of $\pm \infty$).
If the integral converges in a neighborhood of that point then the geodesic
reaches the boundary of $\R$  in a finite
value of the affine parameter and hence the manifold is incomplete: the
projection on the fiber is timelike as follows from the computation of norms.
Then it is easy to see that also this projection reaches the boundary at the
same value of the affine parameter. Indeed let $\phi(s)$ be the
projection of the geodesic $\gamma(s)$ on $\Sigma$ and let $\sigma(s)$ be the
affine (timelike) length of $\phi(s)$ computed with the natural metric of
$\Sigma$; then  
from $-\alpha^2(t(s)) \le(\frac {d\sigma}{ds}\ri)^2=\|V\|^2= C-\|X\|^2$  we
get  
\be
{\rm d}\sigma = \frac {{\rm d}t}{\alpha(t(s)) \sqrt{1+ \frac {\alpha^2(t)
C}{\alpha^2(t_0)(C+ {\dot t_0}^2)} }}\ ,
\ee
which shows that $\sigma$ diverges as $t$ tends towards the infinity where
$\alpha$ vanishes (and this happens at a  finite value of the affine parameter $s$).\\
Similar reasoning holds for the remaining two cases ($C=1,0$); moreover one
realizes that there are no complete inextensible spacelike or lightlike
geodesics.\\
If we change the signature of $\R$ we only have to swap the role of time and
space in the previous derivation.
$\rule{5pt}{8pt}$ \par\vskip 3pt
\bx
In the case $\alpha(t)= \exp(t)$ (which we will encounter later)
 we see that $\alpha$ is integrable at $t=-\infty$. Let $\gamma(s)$ be an
inextensible incomplete geodesic: if it is timelike or null 
then it has no ``turning point'' in $t$ namely $\dot t\neq 0$ and $\dot t<0$
if the geodesic reaches the boundary in the future.\\
If it is spacelike then it has a ``turning point'' and it is incomplete in
both directions.
\ex
\section{Riemannian case}
\label{main}
Throughout this section the manifold $(M,g)$ will be of Riemannian
signature, although Lemmas \ref{firstintegral}, \ref{pregeo},
\ref{totgeo} and much of the proof of Thm. \ref{reduction} 
apply without modifications to the pseudo-Riemannian case
as well.
Given a smooth function $\o:M\to \R$
we will denote the (contravariant) gradient of $\o$ by $\O$.
 The Hessian $\Ho$ will be the second
covariant differential w.r.t. the Levi--Civita connection of $(M,g)$,
 namely 
\be
\Ho(X,Y) =
<\nabla_X\O,Y>.
\ee
\bl
\label{firstintegral}
If $\o$ is an Obata's function
 then 
\be
\|\O\|^2+\k\o^2=\h\label{firstint}
\ee
 for some constant $\h\in \R$.
\el
{\bf Proof}.  Taking derivative along the vector $X$  of $\|\O\|^2+\k\o^2$ we
obtain
\be
X\le(\|\O\|^2+\k\o^2\ri) = 2\Ho(X,\O) + \k\o g(X,\O) = 0
\ee
 This result does not rely on the signature and
holds in the pseudo--Rie\-man\-nian case as well. $\rule{5pt}{8pt}$
 \par\vskip 3pt 
\bl
\label{pregeo}
If $\o$ is an Obata's function then 
the curves generated by $\O$ are pregeodesics.
\el
{\bf Proof}.  Obata's equation $\Ho(X,Y)= g(\nabla_X\O,Y)=-\k \o g(X,Y)$ is
equivalent to $\nabla_X\O = -\k\o X$. In particular
$\nabla_\O\O = -\k \o \O$, and hence $\O$ is parallel. Once more this result does
not rely on the signature and 
holds in the pseudo--Rie\-man\-nian case as well.  $\rule{5pt}{8pt}$ \par\vskip 3pt
\bl
\label{totgeo}
If $\o$ is an Obata's function on a pseudo--Riemannian manifold
$(M,g)$ with first integral $\h\neq 0$ 
and if $\Sigma_0:= \o^{-1}(0)$ is not empty, then it is a smooth
totally geodesic surface.
\el
{\bf Proof}. From the equation 
\be
\|\O\|^2\bigg|_{\Sigma_0} = \h \neq 0
\ee
follows that $\Sigma_0$ is smooth. Let $p\in  \Sigma_0$ and $X\in T_p
\Sigma_0$ and consider the geodesic starting at $p$ with initial
tangent vector $X$, denoted by $\gamma(t)$, $t$ being the affine
parameter. A simple computation gives
\be
\frac {{\rm d}^2}{{\rm d}t^2} \omega (\gamma(t)) =\frac {{\rm d}}{{\rm
d}t} g(\O,\dot\gamma) = g\le(\nabla_{\dot\gamma}\O,\dot\gamma\ri) =
-\varkappa \omega(\gamma(t)) \|\dot \gamma\|^2\ .
\ee
Therefore the function $\chi(t):= \omega(\gamma(t))$ satisfies a
second order ODE $\chi'' = \pm \varkappa \chi$ or $\chi''=0$ according
to the type space, time or light of the geodesic. In our case $\chi(0)=
\chi'(0) = 0$ and hence $\chi(t) \equiv 0$, so that the geodesic
remains in the level surface $\Sigma_0 = \omega^{-1}(0)$.  $\rule{5pt}{8pt}$ \par\vskip 3pt
\bl
\label{maxint}
If $\o$ is an Obata's function  with first integral $\h= \|\O\|^2 +
\k\o^2$ on the geodesically complete Riemannian manifold 
$(M,g)$  then ${\overline J}=\o(M)\subset \R$ is the closure of the interval 
\begin{enumerate}
\item[(i)] $J=\le(-\sqrt{\h/\k},\sqrt{ \h/\k}\ri)$
 if $\k,\h>0$;
\item[(ii)] $J= \le(\sqrt{{|\h|}/\k},\infty\ri )$ (or $J=
\le(-\infty,-\sqrt{{|\h|}/{\k}}\ri)$) if $\h\leq 0$, $\k>0$;
\item[(iii)] $J=\R$ if  $\h>0$, $\k\leq 0$.
\end{enumerate}
\el
{\bf Proof}.
The cases $\h,\k<0$ or $\h=\k=0$ cannot occur in the Riemannian case
because of the positiveness of the metric and eq. (\ref{firstint}).\\
Let $p\in M$ be any point where $\O_p\neq 0$. Consider the geodesic
$\gamma(t)$ starting at $p$ and parallel to $\O$ (by Lemma
\ref{pregeo}), $t$ being the affine length parameter.
 The function $f(t) := \o(\gamma(t))$ then satisfies (from eq. (\ref{firstint}))
\be
(f'(t))^2+\k (f(t))^2 = \h. 
\ee
By virtue of the geodesic completeness,
 this is valid over the whole interval $t\in
\R$. Integrating this simple ODE one obtains
\be
f(t) =\le\{
\begin{tabular}{ll}
$\displaystyle{
\sqrt{\h/\k}\, \cos\le(\sqrt{\k} t\ri) }$
&$ \k>0,\ (\h>0)$, \\
$\displaystyle{
\sqrt{\h/\k}\, \sinh\le(\sqrt{|\k|} t\ri)} $&$ \k<0,\ \h>0
$,  \\
$\displaystyle{
\sqrt{ {|\h|/\k}}\, \cosh\le(\sqrt{|\k|} t\ri)} $&$ \k<0,\ \h<0
$, \\
$\displaystyle{|\k|^{-\frac 1 2}
\exp \le(\sqrt{|\k|} t\ri)} $&$ \k<0,\ \h=0
$,\\
$\displaystyle{
\sqrt{\h}}\,t $&$ \k=0,\ (\h>0)
$.
\end{tabular}\ri.\label{cases}
\ee
The intervals $\overline J$ are just the ranges of the function $f(t)$ in the
 different cases. Note that the open intervals $J$ are constituted
 by all regular values of $\o$. $\rule{5pt}{8pt}$ \par\vskip 3pt
We can now state the theorem in the Riemannian case: the proof was given by
the authors in \cite{kn:bg}. For completeness we report the
derivation of  the result since part of the proof applies without changes to
 the pseudo--Rie\-man\-nian  case as well.
\bt\label{reduction}
Let $(M,g)$ be any complete smooth Riemannian 
manifold of dimension greater
than one (to avoid trivialities) such that there exists an Obata's function
$\o$ with first integral
 $\k\o^2+\|\O\|^2=\h$: denoting by $\Delta:=\{ x\in M\ ;
\langle\O_x,\O_x\rangle =0\}
$ the critical fibers of $\o$,  then\\
i) $(M\setminus \Delta ,g)$ is  isometric to a
warped product $I\times_\alpha\S_\q$ where $I\subseteq \R$ is an open
interval, $\S_\q:=
\o^{-1}(\q)$ for a regular value $\q$, and $\alpha$ is a suitable
function to be specified in the proof.\\
ii) if $\Delta\neq \emptyset$ then $(M,g)$ is of constant curvature
$K^{(M)} =\k$;\\
iii) If $\k\leq 0\leq\h$ then the above holds globally (and $I=\R$).
\et
{\bf Proof}.
First of all, from eq. (\ref{firstint}) and from Lemma \ref{maxint} 
follows that  $\Delta$ is not empty in correspondence of the points
$p\in M$ that solve $
\o^2(x) = \frac \h\k$. This implies that the singular fibers are level
surfaces of $\o$ corresponding to the values $\pm\sqrt{\h/\k}$; since $\o$ is
real--valued, this can occur only if $\h\k\geq 0$.\par
In any case, let now $\q$ be a regular value of $\o$ and $J\subseteq \R$ the
maximal (open) interval of regular values containing $\q$; this is the interval
defined in Lemma \ref{maxint} according to the various cases.\\
The foliation induced on $\o^{-1}(J)\subseteq M$ by the function $\o$
is then regular, all level set of  $\o$
being diffeomorphic by means of the diffeomorphism generated by the
gradient $\O$ of $\o$
\be
 J\times \Sigma_\q : \mathop{\longrightarrow}^{\psi}
\omega^{-1}(J)\subseteq M\ .
\ee
More explicitly, the point $\psi(\tilde q, \sigma)$ is the (unique)
point of $\Sigma_{\tilde \q}:= \o^{-1}(\tilde\q)$ lying on the
geodesic generated by $\O$ and starting at $\sigma\in \Sigma_\q$.
Below we will denote by $\o$ both the function and the coordinate on $J$.
This definition implies (tautologically) that 
$\psi_*\pa_\o= \frac {\O}{\|\O\|^2}$.
We  now prove that the metric  $\tilde g:=\psi^* g$ gives $J\times\S_\q$ the
structure of warped product. 
Let  
 $p_1$ and $p_2$ denote the projections onto the two factors of
$J\times\S_\q$ and   $i:\S_\q\hookrightarrow M$ be the natural
injection. For all
$X,Y$ in the tangent bundle of $\S_\q$ in $J\times\S_\q$ (i.e. in
$\Gamma(p_2^*T \S_\q)$) 
\beas
& &\tilde g(\pa_\o,\pa_\o)=\frac 1{\|\O\|^2} = \frac 1{\h-\k\o^2} 
\ ;\qquad\tilde
g(\pa_\o,X)=0; \qquad
\tilde g(X,Y)=g(\psi_*X,\psi_*Y)\ .
\eeas
Let  $X,Y\in \Gamma(p_2^*T\S_\q)$ such that $[\pa_\o,X]=[\pa_\o,Y]=0$
and thus $[\O,\psi_*X]=[\O,\psi_*Y]=0$:
if we compute ${\cal L}_{\pa_t}\tilde g$ we get
\beaq
& &\pa_\o(\tilde g(X,Y))=({\cal L}_{\pa_\o}\tilde g)(X,Y) =
\frac 1 {\|\O\|^2} \O \bigg(g(\psi_*X,\psi_*Y)\bigg)=\cr
& &=\frac 1 {\|\O\|^2}\le \{g(\D_\O \psi_*X,\psi_*Y)+g(\psi_* X,\D_\O
\psi_* Y)\ri \}=\cr 
& &= \frac 1 {\|\O\|^2}\le\{g(\D_{\psi_*X}\O,\psi_*Y)+g(\psi_*
X,\D_{\psi_* Y}\O)\ri \}=\frac 1 {\|\O\|^2} {2 \Ho(X,Y)}=\cr
& &=\frac{-2\k\o}{\h-\k\o^2} \, g(\psi_*X,\psi_*Y) =  \le(\pa_\o \ln
\le| \h-\k\o^2\ri|\ri)  \,
\tilde g(X, Y) \ ,\label{warped}
\eeaq
That is the metric on $\{\omega\}\times \Sigma_\q$ undergoes conformal
rescaling under change of the base--point $\o$.
Hence 
\bes
\tilde g= \frac 1{\h-\k\o^2}{\rm d}
\o^2+\le|\frac{\h-\k\o^2}{\h-\k\q^2}\ri| \,i^*g
\ees
This proves the warped structure; introducing a geodesic coordinate
$t$ in $J$ according to  
\be
{\rm d}t^2 = \frac 1{\h-\k\o^2}  {\rm d}\o^2
\ee
we obtain the metric $\tilde g$ in the form 
\be
\tilde g= {\rm d}t^2 + \alpha ^2(t) \, i^* g\ ,\qquad
\alpha(t) := \sqrt{\le|\frac{\h-\k\o^2}{\h-\k\q^2}\ri|} = \frac 1
{\sqrt{|\h-\k\q^2|}} |f'(t)|\ , 
\ee
and $f(t)$ being given in eq. (\ref{cases}). According to the different
cases  the intervals $J$ of Lemma \ref{maxint} expressed  in the
geodesic coordinate $t$ are \\ 
(1): $ \k>0,\ (\h>0)$,   $I=(0, \pi/\sqrt{\k})$; \\
(2): $\k<0,\ \h>0$, $I=\R$;\\
(3): $ \k<0,\ \h<0$,  $I=(0,\infty)$;\\
(4): $ \k<0,\ \h=0$,  $I=\R$ ;\\
(5): $ \k=0,\ (\h>0)$,   $I=\R$.\\
We remark that the critical set $\Delta$ is not empty 
only in cases (1) and (3) above, since then the function
$\o$ takes on the critical values as seen in Lemma \ref{maxint}.
Furthermore it follows from elementary Morse theory that since at
those points the Hessian $\Ho_{|_\Delta} = \pm\sqrt{|\h/\k|}
g_{|_\Delta}$  is nondegenerate and definite (positively or
negatively), the critical points are isolated and are either maxima or
minima of $\o$.  
In case (1), i.e.,  $\k>0$, $\h>0$, $\Delta$ is
constituted by two isolated points, 
 one maximum and one minimum $p_\pm$ with the
critical values $\o_{cr}=\pm
\sqrt{\h /\k}$;\\
in case (3), i.e.,  $\k< 0>\h$, $\Delta$ is just one isolated
point of minimum and $\o(M)= [\sqrt{|\h/\k|},\infty)$ (or maximum if
$\o(M)=(-\infty, -\sqrt{|\h/\k|}]$)  with critical value 
$\o_{cr} = \sqrt{|\h|/|\k|}$  ($\o_{cr} = -\sqrt{|\h|/|\k|}$).\\
Moreover, 
since the
Hessian at the critical points is of definite signature (being proportional to the Riemannian
metric), then the level surfaces of $\o$ are topological spheres for
values near to the critical values, as follows once more from
elementary Morse theory. But since all regular level surfaces are
diffeomorphic then all the level sets $\S_\q = \o^{-1}(\q)$ are
topological spheres.\\
Outside of the critical locus $\Delta$ the level surfaces of $\o$
are the same as those  of $\alpha$ (spheres):  as $\o$ tends to a
critical value $\alpha$ tends to zero (by its definition).
To prove assertion (ii)  we
now  compute the sectional curvature of $M$
on a plane spanned by $U, V\in \Gamma(T\S_\q)$. The
calculation follows from  the expression of the 
sectional curvature of a warped product;
\be
K^{(M)}_{UV} = { {K^{\S_\q}_{UV}-(\alpha')^2  }\over{\alpha^2}}\ .
\label{curvv}
\ee
We shrink this topological
sphere 
\bes
\alpha \mathop{\longrightarrow}_{t\to t_0} 0\ \ \Leftrightarrow\ \ \omega\to \omega_{cr}\ ,
\ees
 by parallel translating the two vertical vectors $U,V$ up to the
critical point $p_{cr}$ along the  flow generated by the gradient $\O$ 
(recall that $\O$ and the gradient of $\alpha$ generate pre-geodesics).
For each such  flow line $\gamma$, the projection on the fiber
$\S_\q$ is constant, and the 2-plane spanned by $U,V$ does not change
(each vector is just rescaled). At the end of this shrinking process
we obtain two  vectors in the tangent space
$T_{p_{cr}}M$.
 Since we must obtain a well definite value of the
sectional curvature of $M$ then we must have
$K^{\Sigma_\q}_{UV}=\le(\alpha'(t_0)\ri)^2$ independently of the ``direction''
of the geodesic, namely of the point on $\Sigma_{\q}$, and of the
two-plane. This proves that $\Sigma_\q$ is a sphere because we have
just proved that its sectional curvature is a constant.\\
Then, from Eq. \ref{curvv} and from the explicit form of $\alpha$,
 it follows that the sectional curvature of the manifold $M\setminus \Delta$ is also
constant $K^{(M)}=\k$ and
hence $(M,g)$ is globally (by continuity) of constant sectional
curvature, which proves part (ii). Notice that the fact that $\o$ is a
Morse function with spheres as level sets, implies here  that
$(M,g)$ is actually a round sphere or a non-quotiented hyperboloid
(depending on the sign of $\k$).
\par\vskip 3pt
If the critical locus $\Delta$ is empty (which corresponds to 
the remaining cases 2), 4), 5))  we have no constraint on the
curvature of the leaf $\Sigma_{\q}$, which can be any complete smooth 
Riemannian manifold.
 $\rule{5pt}{8pt}$\par\vskip 5pt
\section{The pseudo--Riemannian case}
\label{mainps}
We now consider Obata's  equation (\ref{magn}) assuming that $g$ is
pseudo--Rie\-man\-nian with signature $(r,p)$ (both non-zero).
We still have the same first integral as in Lemma (\ref{firstintegral}), but
now the square--norm of $\O$ can be of any
sign and hence any combination of signs of $\k$ and $\h$ is a priori allowed.
From eq. \ref{firstint} 
we have the following implications on the type of the vector $\O$ depending on
the relative signs of these two constants.\par\vskip 5pt
\begin{tabular}{lll}
  $\k$  & $\h$ & {\bf Type of} $\O$\\ 
 $\k>0$, & $\h>0$ & depends\\
 $\k>0$, & $\h<0$ & timelike\\
 $\k>0$, & $\h=0$ & timelike or null\\
 $\k<0$, & $\h>0$ & spacelike\\
 $\k<0$, & $\h<0$ & depends\\
 $\k<0$, & $\h=0$ & spacelike or null\\
 $\k=0$, & $\h>0$ & spacelike\\
 $\k=0$, & $\h<0$ & timelike\\
 $\k=0$, & $\h=0$ & null.
\end{tabular}\par
From Lemma \ref{firstintegral} and Eq. (\ref{firstint}) 
it follows that the critical {\em values} of
$\o$ can be only $\pm \sqrt{\frac {|\h|}{|\k|}}$: we denote with $\Sigma_\pm$
the corresponding singular level--sets. These, contrarily to the Riemannian
case, do not coincide with the set of stationary points because the latter are
now saddle--points.\par
In order to extend Thm. (\ref{reduction}) to the pseudo Riemannian
case we will formulate Thms. (\ref{mainpseudo1}, \ref{mainpseudo2}, \ref {mainpseudo3},
\ref {mainpseudo4}) rather than one single theorem in which we have
many subcases.
 Only the case $\k\geq 0$ will be addressed because the
case $\k<0$ can be  obtained easily from the case $\k>0$ by exchanging the
roles of ``time'' and ``space'', that is by swapping $r$ and $p$ in
the signature.
Thus, for instance, the case $\k>0,\ \h>0$ corresponds to the case $\k<0,\
\h<0$ and so on.\\
We will work under the following 
\par\vskip 3pt
{\bf Common assumptions}: 
 {\em 
 $(M,g)$ is a geodesically complete, connected pseudo--Rie\-man\-nian manifold;
the signature of the metric is $(r,p)$ (with $r$ negative and $p$
positive eigenvalues, both nonzero) and 
 there exists a nontrivial Obata's function $\o$
 with first integral $\|\O\|^2+\k\o^2 = \h$.}\vskip 4pt
\bt
\label{mainpseudo1}
Under the  common  assumptions, let $\k>0$ and $\h>0$.
 Let $\Sigma_\pm =
\o^{-1}\le(\pm \sqrt {\h/\k}\ri)$ be the singular fibers.
Then
  $\Sigma_\pm$ are both nonempty and
  $(M,g)$ has constant sectional curvature $K=\k$.\par 
 In each connected component  of 
 $M\setminus \Sigma_+\cup \Sigma_-$ the type of
 $\O$ is constant and $\o$ takes value in one of these intervals 
\be
J_-=(\infty,-\sqrt{\h/\k})\ ,
\qquad  J_0 = (-\sqrt{\h/\k},\sqrt{\h,\k})\ ,\qquad J_+ = (\sqrt{\h/\k},+\infty)\ .
\ee
We denote with $M_{cc}^{+}$,  $M_{cc}^{-}$,  $M_{cc}^{0}$
the generic connected components  of $\o^{-1}(J_+)$, $\o^{-1}(J_-)$,
 $\o^{-1}(J_0)$ .\\ The boundary of $M_{cc}^{\pm}$ is a
 light--cone with vertex in a critical point $p_\pm$ of $\o$, with $\o (p_\pm)
 = \pm \sqrt{h/\k}$. Moreover $M_{cc}^{\pm}$ is the set of points which
 are timelike related to $p_\pm$. 
The boundary of $M_{cc}^{0}$ is constituted by two light--cones with vertex at
 $p_\pm$ and the interior points are spacelike related to $p_\pm$.\\
Each connected component  is 
isometric to an appropriate warped product $I\times_\alpha \Sigma$ as
 indicated below, where $d(\cdot,\cdot):M\times M\to \R$ denotes the
 pseudo--distance (negative for timelike separation):\\
a) in $M^{0}_{cc}$ 
\bea
&& \o(p)= \sqrt{\h/\k}\cos\le(\sqrt{\k} d(p,p_+)\ri)\ ;\qquad I=(0,\pi/\sqrt{\k});\qquad 
\alpha = \sqrt{\h} \sin\le(\sqrt{\k} t\ri)\ ,\\
&& {\rm d} s^2 = {\rm d}t^2
+\alpha^2(t) {\rm d}s^2_\Sigma 
\eea
and $\Sigma$ is a hypersurface of type $(r,p-1)$ and constant curvature
$K^\Sigma= \h$.\\
b) in $M^{\pm}_{cc}$ we have
\bea
&& \o(p)= \pm \sqrt{\h/\k}\cosh\le(\sqrt{\k} d(p,p_\pm)\ri)\ ;
\ I=(0,\infty);\ 
\alpha = \sqrt{\h} \sinh\le(\sqrt{\k} t\ri)\ ,\\
&& {\rm d} s^2 = -{\rm d}t^2
+\alpha^2(t) {\rm d}s^2_\Sigma
\eea
and $\Sigma$ is a hypersurface of type $(r-1,p)$ and constant curvature
$K^\Sigma=-\h$.\\
c) the singular fibers $\Sigma_\pm$ are union of 
light--cones with vertices in the
 stationary points of $\o$.
\et
{\bf Proof}. 
We first prove that $\Sigma_\pm$ are both nonempty.
 Let $p$ be a
regular point of $\o$ where $\O$ is not null; such a point exists otherwise
$\|\O\|^2\equiv 0$ but then eq. (\ref{firstint}) would imply that
$\o$ is  a constant, in contradiction with eq. (\ref{magn}). 
 Let $\gamma(s)$ be
 the geodesic generated by $\O$ starting at $p$.
Then an argument similar to the Riemannian case shows that the value of $\o$
along this geodesic will eventually reach one of the two values $\pm
\sqrt{\k/\h}$. More precisely, if $\O$ at $p$ is spacelike, then
$\o(\gamma(s))$ is a trigonometric function which attains both values
along the geodesic $\gamma$.
 If $\O$ is  timelike at $p$, then $\o(\gamma(s))$ is a hyperbolic
cosine which reaches one of the two values along the geodesic.\\
Let $p_\pm$ be the first points on such a geodesic where
$\o(p_\pm)=\pm\sqrt{\h/\k}$. Clearly $p_\pm\in\Sigma_\pm$ and hence they 
 are not empty.\\
We now prove that $p_\pm$ are critical points. Suppose by contradiction that
$\O(p_\pm)\neq 0$; then it should  be a null vector from
eq. (\ref{firstint}), which is impossible
because the geodesic  $\gamma$ under consideration (generated by $\O$)
is either timelike or spacelike. This argument shows also  that any
point $p$ in $M\setminus \Sigma_+\cup\Sigma_-$ is geodesically
connected with one or both type of the critical points $p_\pm$ (according the the cases). 
Therefore $p_\pm$ are isolated critical points of $\o$ (the Hessian being 
non-degenerate). Moreover, they
are saddle points and thus $\O$ is of every type in any neighborhood
of the critical points since
$\|\O\|^2=\h-\k\o^2$ takes on positive and negative values in the
neighborhood of the critical point.\par 
 It is not
difficult to 
show now that in the connected regions of $M\setminus \Sigma_+\cup
\Sigma_-$ (where $\O$ is spacelike or timelike) 
\bea
&& \o_{space}(p)= \sqrt{ {\h}/{\k}} \cos\le( \sqrt{\k}\,d(p,p_+) \ri)
\label{unno} 
\\
&& \o_{time}(p) = \pm\sqrt{\h/ \k} \cosh\le( \sqrt{\k}\,d(p,p_\pm)
\ri)\ .\label{duue}
\eea
In eq. (\ref{unno}) $p$ is restricted to belong to 
 the region of points which are spacelike  w.r.t. $p_+$ while in
eq. (\ref{duue}) $p$ is timelike related to $p_\pm$.\\
The singular fibers of $\o$ are constituted by the light--cones
through the points $p_\pm$.
 Moreover, following any 
 spacelike geodesic emerging from $p_+$ we reach a stationary point $p_-$ and
vice versa: of course in general 
we need not reach the same stationary point where we
started from, as this depends on the global topology of the manifold which
cannot be fixed in this context (see Remark \ref{nonsimply}).\\
In order to  prove that the  sectional
curvature is  constant we can adapt 
 the argument used in part (ii) of the proof of Thm. \ref{reduction}
 except that now we must approach the stationary points
$p_\pm$ from spacelike or timelike directions.
In these two cases the geometry  of the level surfaces of $\o$ is different,
since they are (locally) modeled on hyperboloids of different signatures.
Nonetheless this is sufficient to prove that the sectional
curvature of the surfaces $\o^{-1}(\q)$ are constants: this in turn forces the sectional
curvature of $M$ to be constant and equal to $\k$ on the complement of the light--cone.
Smoothness of $(M,g)$ completes the proof. $\rule{5pt}{8pt}$ \par\vskip 13pt
\epsfxsize=7cm
\epsfysize=6cm
\centerline{\epsffile{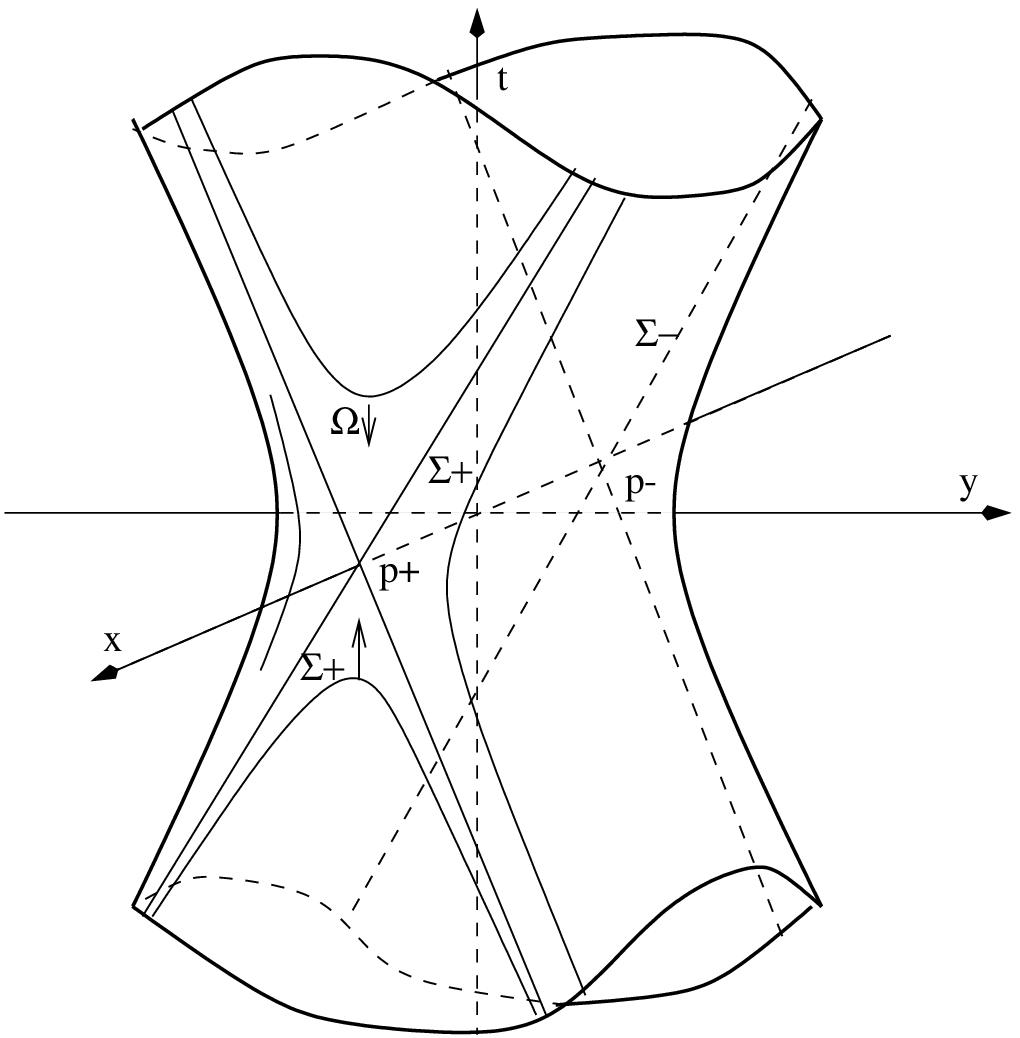}}
\begin{center}
{An example for $k=1$, $h=1$ 
 on the quadric $-t^2+x^2+y^2=1$ with $\o= x$: this is the
paradigmatic case up to a covering.}
\end{center}
\vskip 6pt
\br
\label{nonsimply}
In the case $\k>0$, $\h>0$, (or $\k<0$, $\h<0$) we have proved that 
$(M,g)$ must be of constant sectional curvature and now  we can
describe $\o$ explicitly.
If we realize $(M,g)$ as a suitable quadric (or a
covering of it), then $\o$ is just any linear function in the ambient
pseudo--Rie\-man\-nian
 flat manifold restricted to the quadric, provided that its
gradient is not null.\\
In this case we have more than two critical points
only if the manifold is a covering of the above quadric. 
Consider for instance 
the quadric $-{Z}^2-{Y}^2 + \sum_1^d {X^j}=-1$ in  a flat spacetime with
metric $-{\rm d}Z^2 - {\rm d}Y^2 + \sum_1^d {\rm d}{X^j}^2$ 
 (called ``Anti de Sitter spacetime'') and the function
$\o = Z$, (this corresponds to the case $\k=-1,\  \h=-1$): since this space
has a nontrivial fundamental group $\pi_1 \simeq{\mathbb Z}$,
 we may pass to its
universal covering (or other coverings), and the function $\o$ would have many
critical points.
\er
We now consider the case $\k>0$ and $\h<0$. Contrarily to the previous
one
we do not obtain a complete rigidification and we are left with an 
arbitrariness in the metrics that may occur on the leaves of Obata's function.
This is similar to what happens in the Riemannian case with
$\k<0,\ \h>0$.\\[10pt]
\epsfxsize=7cm
\epsfysize=6cm
\centerline{\epsffile{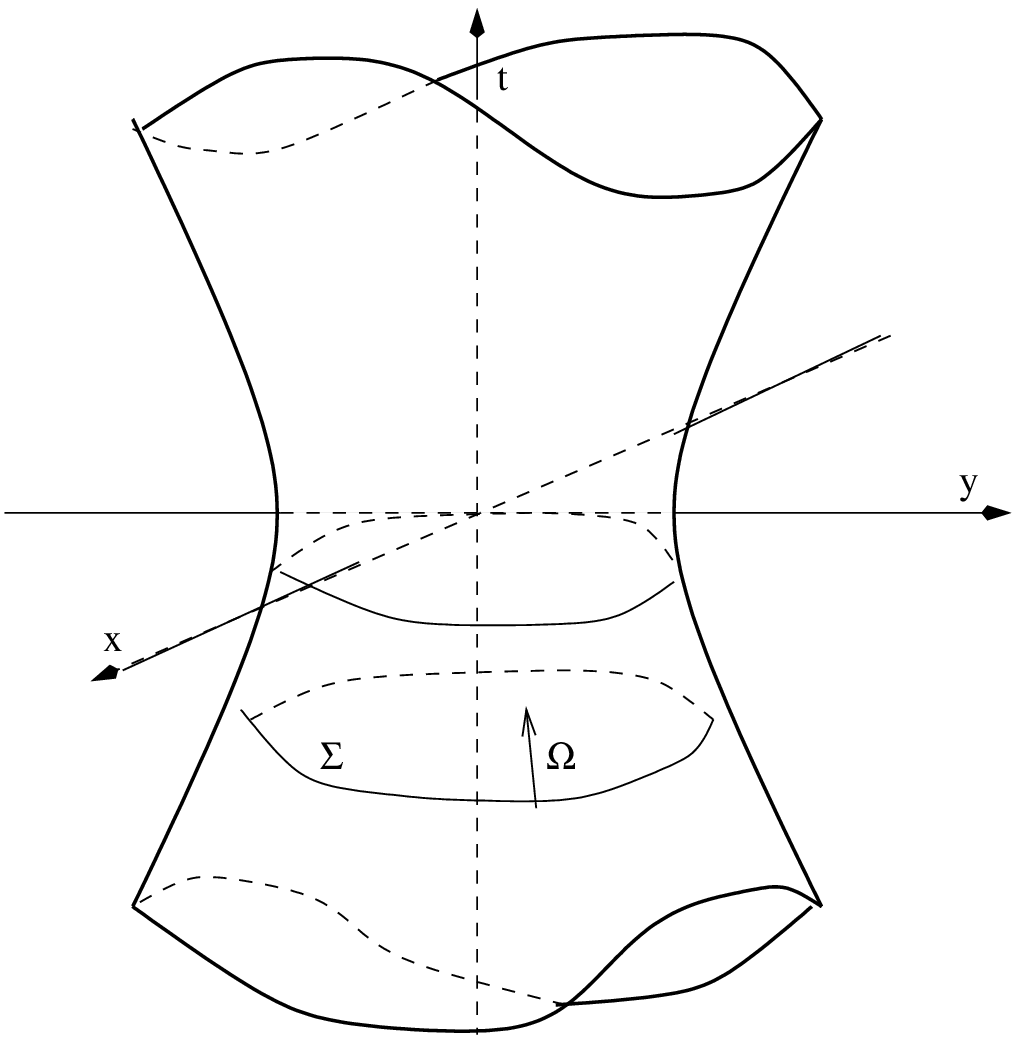}}
\begin{center}
\baselineskip 14pt plus 2pt minus 2pt
An example for $\k=1$ and $\h=-1$ 
on the quadric $-t^2+x^2+y^2=1$ with $\o= t$: this is {\bf
not} the paradigmatic case as in general the manifold does not have
constant curvature.
\end{center}
\vskip 5pt  
\bt
\label{mainpseudo2}
Under the common  assumptions, let now  $\k>0$ and $\h<0$.
Then $\O$ is everywhere timelike and $(M,g)$ 
is globally isometric to a warped product $\R\times_\alpha \Sigma$ with
\beaq
&& \o(x) = \sqrt{|\h|/\k}\sinh\le(\sqrt{\k}\,t\ri);\ 
\alpha (t)= \sqrt{|\h|}\cosh\le(\sqrt{\k}\,t\ri)\ ,\cr
&&  {\rm d}s^2 = -{\rm d}t^2 + \alpha^2(t){\rm d}s^2_\Sigma\cr
&& x:= (t,\sigma)\ \in \ \R\times \Sigma
\eeaq
where $\Sigma$ is any geodesically complete, connected 
pseudo--Rie\-man\-nian submanifold of type $(r-1,p)$.
\et
{\bf Proof}. 
In this case  $\O$ is everywhere timelike because $\|\O\|^2 = -|\h|
-\k\o^2\leq -|\h|<0$; in particular this implies that the foliation
induced by the level surfaces of $\omega$ is smooth and all leaves are
diffeomorphic by means of the flow induced by the gradient $\O$. 
By integrating
$\omega$ along the geodesic generated by its gradient $\O$ we find
that it can be written, in term of a suitably shifted affine
parameter, as  
\be
\o(\gamma(t)) =  \sqrt{\frac {|\h|}{\k}} \sinh\le(\sqrt {\k} t\ri)\ .
\ee
In particular $\Sigma_0 := 
\o^{-1}(0)$ is non-empty and hence, by Lemma \ref{totgeo}, it is a
totally geodesic smooth hypersurface of type $(r-1,p)$.
This proves that $M$ is diffeomorphic to $\R\times \Sigma_0$, where we
 use the coordinate $t:=\frac 1{\sqrt{\varkappa}}
 \sinh^{-1}\le(\sqrt{\frac{\varkappa}{|\h|}}\omega \ri)$ for the factor
 $\R$. The same computation as in eq. (\ref{warped}) proves that the
 metric structure is isometric to the warped product 
 $\R \times_\alpha \Sigma$ with $\alpha(t) =
\sqrt{|\h|} \cosh\le(\sqrt{\k}t\ri)$.
Finally, Lemma (\ref{completeness}) adapted to the present situation with $B=\R$
and $F=\Sigma_0$,  proves that $\R\times_\alpha \Sigma_0$ is
geodesically complete and hence so is $(M,g)$. $\rule{5pt}{8pt}$ \par\vskip 3pt
The next case ($\k>0$, $\h=0$) is more complicated: the equation does not
rigidify completely the manifold. The situation is  intermediate
between the complete rigidification which we obtain in
Thm. (\ref{mainpseudo1}) where the level surfaces of $\o$ are forced to have
constant curvature  and the very weak one we have in
Thm. (\ref{mainpseudo2}), where the level surfaces of $\o$ may have any
pseudo--Rie\-man\-nian geometry.
In the present case  geodesic completeness of
$M$ implies a condition of asymptotic flatness along any spacelike direction
in  the level surfaces of $\o$.\\
\bt
\label{mainpseudo3}
Under the common  assumptions, let now  $\k>0$ and $\h=0$.\\
(i) Then $\O$ is almost everywhere 
 timelike. \\
(ii) The two subsets $M_\pm = \o^{-1} (\R^\pm)\subset M$
are   isometric to a warped product $\R\times_\alpha \Sigma_\pm$ with
\bea
&&\o(p_\pm) = \pm \sqrt{1/\k} \exp \le(\sqrt{\k} t\ri);\ 
\alpha (t)= \exp\le(\sqrt{\k}\,t\ri)\ ,\cr 
&&{\rm d}s^2 = -{\rm d}t^2+\alpha^2(t) {\rm d}s^2_{\Sigma_\pm}\cr
&& p_\pm = (t,\sigma_\pm)\ \in\ \R\times \Sigma_{\pm}\ , 
\eea
and $\Sigma_\pm$ are  geodesically complete 
pseudo--Rie\-man\-nian manifolds of type $(r-1,p)$.\\
(iii) Necessary condition for
geodesic completeness is that the sectional curvatures of $\Sigma_\pm$
vanish at least as  ${\cal O}(\sigma^{-1})$  along any spacelike geodesic,
$\sigma$ being the natural length in $\Sigma_\pm$.\\
(iv) If the sectional curvature of $(M,g)$ is bounded,
 then $\Sigma_\pm$ are both
flat and hence the sectional curvature of $M$ is actually constant $K=\k$.
\et
{\bf Proof}. 
 We assume that $\o$ is not
identically zero, hence either one or the other of $M_{\pm} :=
\o^{(-1)}(\R_\pm\setminus\{0\} )$ is not empty: without loss of  generality we
assume $M_+$ not empty.\\
The same reasoning as in the Riemannian case shows that $\o$ has the form in
the statement of this theorem, where for now $\Sigma_+$ is just a geodesically
complete pseudo--Rie\-man\-nian hypersurface of the appropriate type.
Then $( M_+,g)$ is isometric to $\R\times_\alpha
\Sigma_+$, with $\alpha(t) = \exp(\sqrt{\k} t)$ and metric
${\rm d}s^2 = -{\rm d}t^2 + \alpha^2(t){\rm d} s^2_{\Sigma_+}$.
Such a warping function falls into the class of functions considered in Lemma
(\ref{incompleteness}), hence $( M_+, g)$ is not geodesically complete and
 $M_+$ is a proper subset of $M$.\par
We now prove that $M_-$ is nonempty as well and that $\O$ never vanishes so
that each connected component of 
 $\Sigma_0:=\o^{-1}(0)$ is a smooth null hypersurface.\\
Indeed, let $p$ be any
point of $\Sigma_0$ and consider any geodesic $\gamma(s)$ starting from
inside $M_+$ and arriving at $p$ at $s=s_0$.
 Then it is possible to compute explicitly
$\omega(\gamma(s))$ for $s<s_0$ using Eq. (\ref{affineparameter}) and prove that 
\be
{\rm d}\o (p)(\dot \gamma(s_0)) 
= \lim_{s\to s_0^-}  \frac d{ds} \omega(\gamma(s))\neq 0\ .
\ee
%For example, for a geodesic of type light we have $
%\omega(\gamma(s)) =(\sqrt{|\k|} |\dot t_0|)^{-1}\le(s_0-s\ri)$, 
This implies that:\\
i) ${\rm d}\o$ does not vanish at any $p\in \Sigma_0 = \omega^{-1}(0)$
and thus \\
ii) $M_-$ is nonempty.\\
Therefore  $M=M_+\cup \Sigma_0\cup M_-$,
$\Sigma_0$ being a smooth light--like 
hypersurface; smoothness is guaranteed by the nonvanishing of ${\rm d}\o$.\\
The other half $M_-$ is also isometric to a similar warped product $M_-\simeq
\R\times_\alpha \Sigma_-$ with the same $\alpha(t)$: the boundary $\Sigma_0$ 
is at $t=-\infty$ for both.
We now prove that $\Sigma_\pm$ must be asymptotically flat in the spacelike
direction:
this  follows from the requirement that the sectional curvature
of $M$ does not blow up as we follow a geodesic which crosses the boundary of
$M_+$ and from the fact that any
geodesic $\gamma(s)$ which crosses the boundary projects onto a spacelike
pregeodesic $\phi(s)$ 
in $\Sigma_\pm$.\\
Indeed, from Lemma (\ref{incompleteness}) we see that in the case
$\alpha(x(s))= \exp( \sqrt{\k}\,x(s))={\cal O}(s-s_0)$ (when the geodesic is
unbounded towards $t= -\infty$)
 which implies that the length of the pregeodesic
$\phi(s)$ grows as $1/(s-s_0)^2$.\\
Considering the expression of the sectional curvature of
$M_\pm=\R\times_\alpha\Sigma_\pm$ \footnote
{Notice the sign $+$ in front of
$\alpha'$ which  comes from the negative signature of $-{\rm d}t^2$}
\be
K_{UU'}^{M_\pm} = \frac {K^{\Sigma_\pm}_{UU'}+
(\alpha'(t))^2}{\alpha^2(t)} =\frac {K^{\Sigma_\pm}_{UU'}+
\varkappa {\rm e}^{2\sqrt{\varkappa}t} }{{\rm e}^{2\sqrt{\varkappa}t}} = 
K^\Sigma_{UU'}{\cal O}\le( \frac 1{(s-s_0)^2}\ri) +\k \ ,\label{stca}
\ee
we see that $K^{\Sigma_\pm}_{UU'}$,
($U,U'\,\in\,T_{\phi(s)}\Sigma_{\pm}$),  must be infinitesimal
w.r.t. $(s-s_0)^2 = \mathcal O(\sigma^{-1})$ (now $\sigma$ is
 the length parameter of the
pregeodesic $\phi$ in $\Sigma_\pm$),
 namely
\be
K^{\Sigma_\pm}_{UU'}  = {\cal O}\le(\frac 1 \sigma\ri)
\ee
along the pregeodesic $\phi(s)\subset \Sigma_\pm$.\par
Finally, we prove that boundedness of the total sectional curvature
$K$ implies that it is constant.
Indeed, using again eq. (\ref{stca}) on the horizontal geodesic
generated by $\O$ (which projects to a constant in $\Sigma_\pm$) and
then sending $t\to -\infty$, we see that $K^{\Sigma_\pm}_{UU'}$ must
identically vanish in $\Sigma_\pm$ and hence they are flat manifolds.
 In this case then, the sectional curvature of $(M,g)$ is
constant and  equal to $\k$. $\rule{5pt}{8pt}$ \par\vskip 3pt
In the next proposition 
we seek to show that there exist examples in which the
manifolds $\Sigma_\pm$ are not flat but only asymptotically flat
in the spatial directions.\\[10pt]
\epsfxsize=7cm
\epsfysize=6cm
\centerline{\epsffile{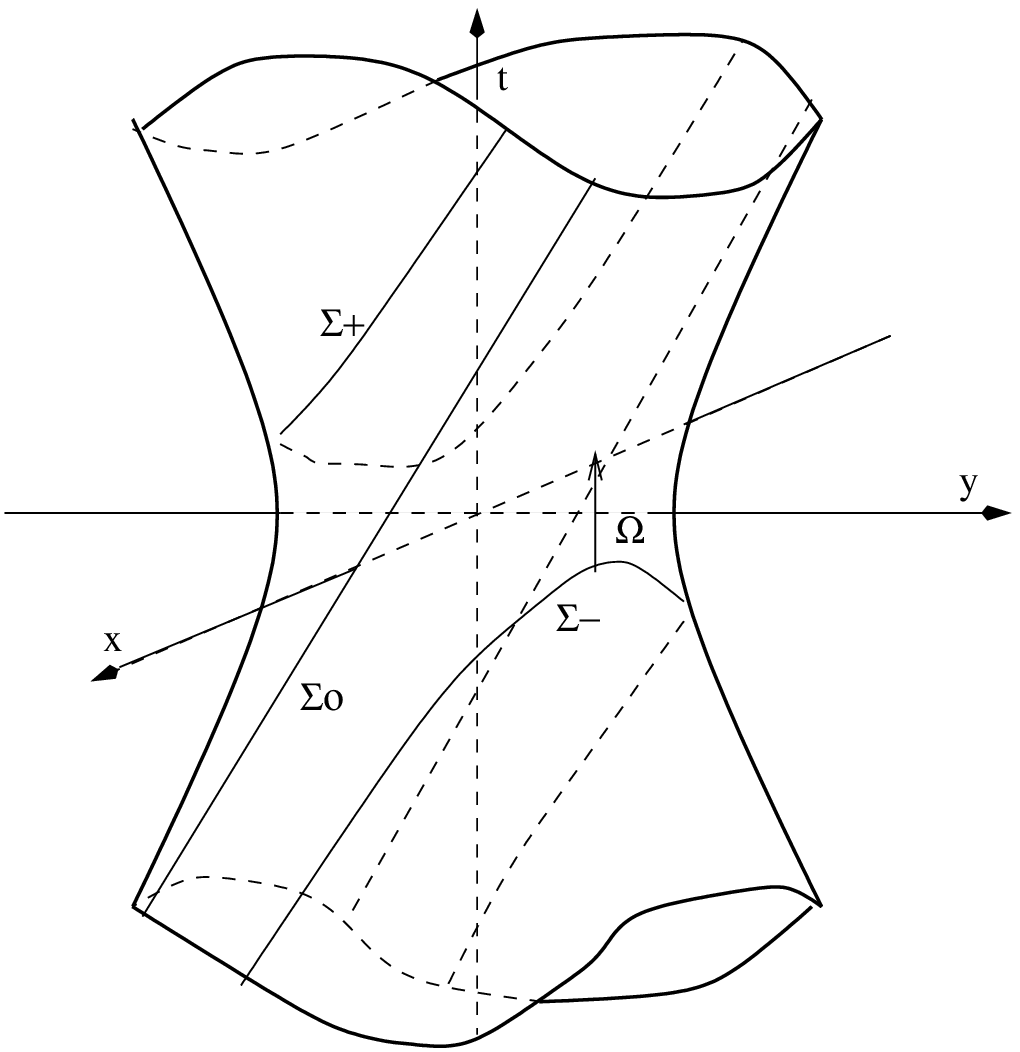}}
\begin{center}
\baselineskip 14pt plus 2pt minus 2pt
An example for $k=1$, $h=0$ on the quadric $-t^2+x^2+y^2=1$ with $\o= t-y$: this is
 {\bf not} the  
paradigmatic case as in general the leaves $\Sigma_\pm$ must only be
asymptotically flat.
\end{center}
\vskip 5pt
\bp
\label{counterex}
Let $\Sigma_\pm$ be two geodesically connected, semi--Rie\-man\-nian
 manifolds of
type $(r-1,p)$ diffeomorphic to  $\R^f$ and such that along any spacelike
geodesic the sectional curvature vanishes faster than any power of $\frac 1
\sigma$, $\sigma$ being the geodesic distance. Set $\alpha(t)= \exp(\sqrt{\k}
t)$ and form the two warped products $(M_\pm, g_\pm)$ as in
\be
{\rm d}s^2_\pm := -{\rm d}t^2 + \alpha^2(t) {\rm d}\sigma^2_{\Sigma_\pm}\ .
\ee
Then we can smoothly glue them in a geodesically complete manifold $(M,g)$
along a null hypersurface $\Sigma_0$: this hypersurface is
connected except in the case of signature $(r,1)$, where it is
 constituted of two connected components.\\
On this manifold the function defined by $\o(t,\sigma) \big |_{M_\pm} =
\pm\frac 1{\sqrt{\k}} \alpha(t)$,
 $\o\big|_{\Sigma_0} \equiv 0$ is smooth and
satisfies $\Ho=-\k\o g$, $\|\O\|^2 = -\k\o^2$.
\ep  
{\bf Proof}. 
We define  the geodesic boundary of $M_\pm$ and extend there  the
metric by means of
the exponential map.
We work on $M_+$, the arguments being identical for $M_-$.\\
Fix an arbitrary point $p_0\equiv 
(t_0,\sigma_0)\in M_+$ and consider the
exponential map $\Exp$ at this point. 
 We saw in Lemma (\ref{incompleteness}) that  there
exist inextensible incomplete geodesics $\gamma(s)$ 
 of any type in $M_+$ starting from $p_0\equiv (t_0,\sigma_0)$ and
that in all these cases the projection of $\gamma$ on $\Sigma_\pm$ is
a spacelike pregeodesic.
Therefore $\Exp:D_{p_0}
\subset T_{p_0}M_+\to M_+$ is defined on a proper star-shaped
subset $D_{p_0}$ of the tangent space. In view of the discussion of the
properties 
of incomplete geodesics, the boundary $\pa D_{p_0}$  of $D_{p_0}$
 is made of three pieces,
according to the type of the incomplete geodesic;
writing $ T_{p_0}M_+\ni X=\pi_*X+\sigma_*X = v+V$, we have the boundary
corresponding to 
the incomplete {\em spacelike}  geodesics 
\beaq
{\mathfrak S}:=\{
X=(v,V)\in T_{p_0} M_+\simeq \R\oplus T_{\sigma_0}\Sigma_+:\cr
 -v^2   + \|V\|^2 ={d_s}^2(X)>0,\ v\neq 0 \}\ ,
\eeaq 
to the incomplete {\em timelike} geodesics
\beaq
{\mathfrak T}:= 
 \{X=(v,V)\in T_{p_0} M_+\simeq \R\oplus T_{\sigma_0}\Sigma_+;\ \|V\|^2>0,\cr
-v^2+\|V\|^2 =-{d_t}^2(X)<0,\  v<0\}\ ,
\eeaq 
and to the incomplete {\em lightlike} geodesics
\beaq
{\mathfrak L}:=
\{X=(v,V)
\in T_{p_0} M_+\simeq \R\oplus T_{\sigma_0}\Sigma_+;\cr
 \|V\|^2={d_l}^2(X)>0,\
 v^2={d_l}^2(X)\}\ .
\eeaq
In these formulas the three functions $d_s,\ d_t,\  d_l$ are the upper 
extrema of the maximal interval of definition of the geodesics and
will explicitly computed in equations (\ref{dl}, \ref{dt}, \ref{ds}). By their
definition $d_s$ and $d_t$ are homogeneous invariant functions of their
arguments because they depend only on the ``direction'' $X/\sqrt{|<X,X>|}$. 
 On the contrary  $d_l$ is homogeneous of degree $-1$ because a null 
geodesic with initial tangent vector $\lambda X$, ($\lambda \in R_+$
and $\|X\|^2=0$) is defined on the interval
$(0,\frac 1 \l d_l(X))$.\\
Consider the set--theoretical union  $\overline M_{+,p_0}:=
M_+\cup \pa D_{p_0}$; we presently 
define a smooth structure of manifold--with--boundary on it.
 Indeed let  $U$ be a  generic open neighborhood of a point $p$ belonging to 
 $\pa D_{p_0}\subset T_{p_0}M_+$:
the corresponding  open neighborhood of $p$ in $\overline M_{+,p_0}$
 will be $(U\cap\pa D_{p_0})\cup \Exp (U\cap D_{p_0})$.
Since  $\Exp$ restricted to
$U\cap D_{p_0}$ is a local diffeomorphism, we have thus   defined
 a smooth structure of
manifold--with--boundary on the set $M_+\cup \pa D_{p_0}$.\\[10pt]
\epsfxsize=6cm
\epsfysize=6cm
\centerline{\epsffile{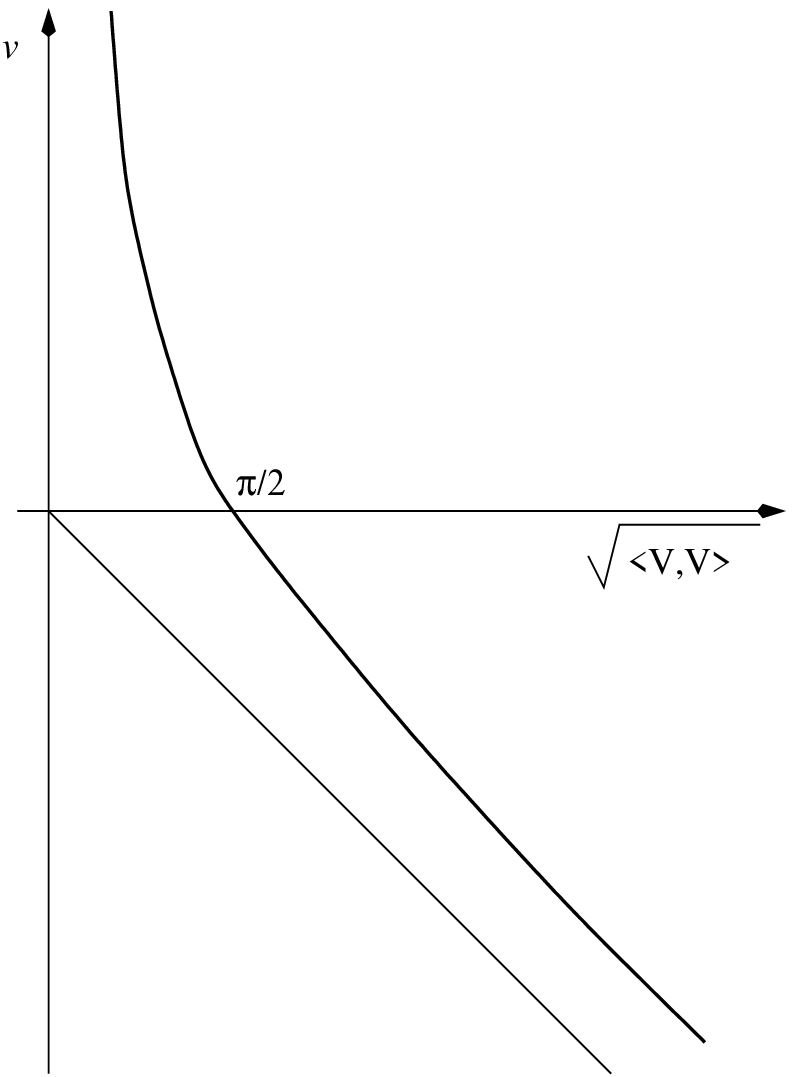}}
\begin{center}
\baselineskip 14pt plus 2pt minus 2pt
The boundary of $M_+$ as appears  via the exponential map
in $T_{p_0}M_+$ here depicted for the special value $\k=1$:
 in abscissa appears the square root of the projection 
$V=\sigma_*X$, which is spacelike for all incomplete geodesics. The curve is
obtained solving the implicit equations 
$\|X\|^2 = -{d_t}^2(X)$ and $\|X\|^2 = {d_s}^2(X)$, where $d_{s,t}$ are given
in the text.
\end{center}\vskip 5pt   
The definition of the boundary seems to rely on the choice
of the base--point $p_0$; we now prove that this is not the case.
Indeed we only need to show that the three functions $d_t,\ d_l,\ d_s$ do not
depend on the point: but this is obvious because, integrating
eq. (\ref{affineparameter})  with $\alpha(t)=\exp(\sqrt{\k}t)$,
 we obtain
\bea
&& d_l(X) = \frac 1 {\alpha (t_0)v}\int_{-\infty}^{t_0}
\alpha(t) {\rm d}t= \frac 1{\sqrt{\k} v}\label{dl} \\
&& d_t(X) = \int_{-\infty}^{t_0}
\frac {\alpha(t) {\rm d}t}{\sqrt{ \alpha^2 (t_0)({\dot t_0}^2- 1) +
\alpha^2(t)}}= \frac 1{\sqrt{\k}} \ln\le(\sqrt{\frac {\dot t_0+1}{\dot t_0-1}} 
\ri)\label{dt}\\
&& d_s (X)= \int_{-\infty}^{t_0}
\frac {\alpha(t) {\rm d}t}{\sqrt{ \alpha^2 (t_0)({\dot t_0}^2+ 1) -
\alpha^2(t)}} =\frac {1}{\sqrt{\k}} \arcsin \le(\frac 1{\sqrt{{\dot t_0}^2+1}}
\ri)\ ,\label{ds}
\eea
where $\dot t_0 = \sqrt{\frac {v^2}{|-v^2+\|V\|^2|}}$.
We can thus remove the subscript $p_0$ and set
 $\overline M_+:= \pa D\cup M_+=: \Sigma_0\cup M_+$, where
$\Sigma_0$ is $\pa D$ when thought in $\overline M_+$.\\
As for the topology of the boundary,
 if the set $\|V\|^2>0$ is disconnected, then so is
$\Sigma_0$; this happens only if the signature of the metric is $(r,1)$
 in which case there are two connected components.\\
We perform a similar construction in $M_-$ and identify the two
copies of the boundary. We therefore obtain a smooth manifold $M=M_-\cup
\Sigma_0 \cup M_+$.\\
We must now define the metric on $TM\big |_{\Sigma_0}$: to this purpose we 
consider the canonical realization of the space of constant curvature $\k$ as
a quadric in a suitable semi-Rie\-man\-nian flat manifold.
Then we can identify $\Sigma_0$ with the intersection of this quadric with a
null hyperplane, and {\em define} the metric on $\Sigma_0$ as the pull-back of
the metric on the quadric.
Since the sectional curvature on $M$ (which is insofar defined only on the
complement of $\Sigma_0$) tends to a constant on $\Sigma_0$ 
with all its derivatives tending
to zero (this follows from the assumptions on the asymptotic flatness of
$\Sigma_\pm$), we have that the metric is smoothly defined also on
$\Sigma_0$.\\
To conclude the proof we only have to show that $\o$ satisfies Obata's
equation: but clearly it satisfies the equation on $M_\pm$ and then by continuity
on the whole $M$. $\rule{5pt}{8pt}$ \par\vskip 3pt
We conclude with
\bt
\label{mainpseudo4}
Under the common  assumptions, let now  $\k=0$ then 
\begin{enumerate}
\item[i)]  if  $h>0$ then $\O$ is everywhere spacelike  and $(M,g)$
is globally isometric to a direct product  $\R\times \Sigma$ with
\bea
\o(x) =  \sqrt{\h}\,t\ ;\qquad 
 {\rm d}s^2 = {\rm d}t^2 +\h {\rm d}s^2_\Sigma\cr
x = (t,\sigma)\ \in\ \R\times \Sigma
\eea
and $\Sigma$ is any geodesically complete 
pseudo--Rie\-man\-nian manifold of type $(r,p-1)$.
\item[ii)]  if $h<0$ then $\O$ is everywhere timelike  and $(M,g)$
is globally isometric to a direct product $\R\times \Sigma$ with
\bea
\o(x) =  \sqrt{|\h|}\,t\ ;\qquad 
 {\rm d}s^2 = -{\rm d}t^2 +|\h| {\rm d}s^2_\Sigma\cr
x = (t,\sigma)\ \in\ \R\times \Sigma
\eea
and $\Sigma$ is any geodesically complete 
pseudo--Rie\-man\-nian manifold of type $(r-1,p)$.
\item[iii)] if $\h=0$ then, apart from the constant solution, $\O$ is a null
Killing vector of $(M,g)$ and the level surfaces of $\o$ are all
totally geodesic.
\end{enumerate}
\et
{\bf Proof}.\\
{\bf Cases i), ii)}.
Each level surface is non--singular: $(M,g)$ is globally
isometric to $\R\times\Sigma$  (with the appropriate types) as a direct
product ($\alpha$ is now  a constant).
 Such a manifold is clearly geodesically
complete iff $\Sigma$ is complete: no further requirement is needed on
$\Sigma_0$ coming from geodesic completeness or smoothness.\\
{\bf Case iii)}.
There is the obvious solution $\o=const$ which is trivial and implies
no requirements whatsoever on the pseudo--Riemannian structure of
$(M,g)$.\\ Let us consider a non-constant solution.
The equation $\Ho\equiv 0\leftrightarrow \nabla\O=0$ proves that $\O$ is
parallel and hence it is  a never vanishing Killing vector.\\
Moreover each level surface is totally geodesic: indeed, if $\varphi(s)$ is a
geodesic starting at a point of $\Sigma_\q$ with initial velocity tangent to
$\Sigma_\q$ then we have
\be
\frac {{\rm d}^2}{{\rm d}s^2} \o(\varphi(s)) = \frac {{\rm d}}{{\rm d}s}
<\O_{\varphi(s)},\dot \varphi(s) = <\nabla_{\dot \varphi}\O,\dot \varphi> +
<\O,\nabla_{\dot \varphi}{\dot \varphi}> = 0\ .
\ee
Hence $\o(\varphi(s)) = \o(\varphi(0)) + A\,s$, but $A$ must vanish because
$A= <\O,\dot\varphi(0)>=0$. $\rule{5pt}{8pt}$ \par\vskip 3pt
A typical example of Case (iii) of Thm. \ref{mainpseudo4}
 is the following: take $\R^{n+2}$ with coordinates
$\o,\xi,\eta_1,..,\eta_n$ and consider the metric 
\be
ds^2 = A\, {\rm d}\o^2 +{\rm  d}\o{\rm d}\xi + {\rm d}\o \sum_{j\geq 1} B_j
{\rm d} \eta_j + \sum_{i,j\geq 1} C_{ij} {\rm d}\eta_i{\rm d}\eta_j\ ,
\ee
where $A, B_i, C_{ij}$ are arbitrary functions independent of $\xi$
(with the only requirement that the metric be nondegenerate).
Then $\O= \frac \pa{\pa\xi}$ and it is a null vector; the condition that $\O$
 is a
Killing vector is ensured by the fact that the coefficients do not
depend on $\xi$.\\
Thus there are not stringent rigidification in this case as well.
\section{Foliations associated to Obata's equation}
\label{maximal}
We now study the implications of the existence of more than one solution: the metric
will be pseudo-Rie\-man\-nian unless otherwise stated.\\
Note that for any Obata's function  $\o$  we have $K_{\O,X} =
\k$. Moreover there cannot exist two Obata's functions corresponding to
different values of the constant $\k$ as we show in the next simple 
\bl
If $\o_1,\o_2$ satisfy $H^{\o_i} = -\k_i\, \o_i g $, $i=1,2$ for some
constants $\k_i$, then $k_1=k_2$.
\el
{\bf Proof}. 
This follows from $K_{\O_1\O_2} = \k_1 = \k_2$. $\rule{5pt}{8pt}$ \par\vskip 3pt
We introduce the following natural
\bd 
 A {\bf maximal system} of Obata's functions is a $m$-tuple $\{
\o_1,...,\o_m\}$ of solutions such that their gradients are almost
everywhere linearly independent.
\ed
The integer $m$ is a pseudo--Rie\-man\-nian invariant of the manifold $(M,g)$.
As we saw in the proof of Thms. \ref{mainpseudo1}, \ref{mainpseudo2},
\ref{mainpseudo3}, \ref{mainpseudo4},if a Obata's function 
 has any critical point, then the manifold must be of
constant sectional curvature: this happens whenever $\k\h>0$ and under
hypotheses of geodesic completeness of the manifold.
Then we can show  that there
exist $n+1$ nontrivial solutions of the equation which 
 are obtained by realizing the manifold in the canonical way as a suitable
quadric in a flat semi-Rie\-man\-nian manifold $\R^{n+1}$ and hence using the
linear 
coordinates of $\R^{n+1}$:
 one of these solution is functionally but not linearly
dependent on the remaining, which then constitute a maximal system in the sense
above.\\
In view of this remark we have the
\bt
The manifold $(M,g)$ has constant sectional curvature if and only if there
exists a maximal system of solutions with $m=\hbox{dim}(M)$.
\et
We now study the intermediate cases $0<m<\hbox{dim}(M)$; if we assume geodesic
completeness then we are addressing only the cases $\k\h_i\leq 0$, $\forall
i=1..m$, otherwise the analysis will only be local.\\
We promptly have 
\bp
The manifold $(M,g)$ is foliated by submanifolds $S^m$  of dimension $m$
with constant sectional curvature $\k$: these foliations are totally geodesic.
\ep
{\bf Proof}.  We use Frobenius's Theorem, showing that $\O_i$ form an
involutive distribution. Indeed 
\be
\le[\O_i,\O_j\ri] = \nabla_{\O_i}\O_j-\nabla_{\O_j}\O_i = -\k\o_j\O_i+
\k\o_i\O_j\ .
\ee
If we have  a geodesic with initial vector $H(p_0)=
\sum_{i=1}^m h^i\O_i(p_0)$ then clearly the vector field $H= \sum_{i=1}^m
h^i\O_i$ generates this geodesic and hence the distribution is also totally
geodesic: this proves that the intrinsic sectional curvature $K^{S^m}$ equals
the sectional curvature of $M$ and hence 
\be
K^{S^m}_{\O_i\O_j} = K^M_{\O_i\O_j} = \k\ ,
\ee
which ends the proof. $\rule{5pt}{8pt}$ \par\vskip 3pt
Before studying the foliation we anticipate the 
\bl
If the dimension $m$ of a maximal system 
 is strictly less than the dimension of
$M$, then any other solution is a linear combination of the basis
$\o_1,...,\o_m$. 
\el
{\bf Proof}. 
Let $\o$ be another solution, then from the assumption of maximality  $\O=
a^i\O_i$ for some functions $a^i$: we are to prove that $a^i$ must be
constants.
Indeed
\be
 -\k\o X = \nabla_X\O = <X,A^i> \O_i +a^i\nabla_X \O_i = <X,A^i> \O_i  -\k
a^i\o_i X\ .
\ee
This implies that 
\be
\k(a^i\o_i-\o) X = <X,A^i> \O_i\ .
\ee
Taking an arbitrary vector field $X$ not belonging to the span of the $\O_i$'s
we must have $\o\equiv a^i\o_i$, so that now $<X,A^i>\O_i\equiv 0$ and hence,
being $\O_i$ independent, $A^i\equiv 0$, namely $a^i$ are constants. 
$\rule{5pt}{8pt}$ \par\vskip 3pt
This lemma shows that if the maximal system is not total (i.e., $m<n$) then we lose
one solution which should be functionally but not linearly dependent on the
other. \\
We remark that 
\bp
The following formula holds
\be
<\O_i,\O_j> = -\k \o_i\o_j+ c_{ij}\ ,
\ee
for some constants $c_{ij}$.
\ep
{\bf Proof}. 
We have
\beaq
&& X<\O_i,\O_j> = <\nabla_X\O_i,\O_j> +<\nabla_X\O_j,\O_i> =\cr 
&& =-\k\o_i<X,\O_j>
-\k \o_j<X,\O_i> = -k X(\o_i\o_j)\ ,  
\eeaq
hence $<\O_i,\O_j> +\k \o_i\o_j=c_{ij}$ is constant. \\
Moreover we can assume that the matrix $c_{ij}$ is diagonal up to a
linear change of basis $\tilde \o_i$. $\rule{5pt}{8pt}$ \par\vskip 3pt
We finally introduce the complementary foliation $\cal F$ to the distribution
spanned by $\{\O_i\}_{i=1..m}$. Clearly the fibers $\Sigma_{p_0}$  of $\cal F$ are the joint
level sets of $\{\o_1,...,\o_m\}$, namely $\Sigma_{p_0} = \bigcap_{i=1}^m
\o_i^{-1}(\o(p_0)) =  \bigcap_{i=1}^m \Sigma_{p_0}^i$
 whose second fundamental form is given by 
\be
\alpha(X,Y) = -\sum_{i=1}^m \frac {H^{\o_i}(X,Y)}{\|\O_i\|^2} \O_i =-\k <X,Y>
\sum_{i=1}^m \le(\frac {\o_i}{\|\O_i\|^2} \O_i\ri) \ ,\label{cucco}
\ee 
for any $X,Y\in T\Sigma_{p_0}$. Eq. (\ref{cucco}) that $\Sigma_{p_0}$ is totally
umbilical.\\
We now restrict to the Riemannian case and, referring to the classical 
definitions in \cite{kn:mo}, we have
\bt
$\cal F$ is a Riemannian foliation and $g$ is bundle--like.
\et
{\bf Proof}. 
Consider the connection $\tilde \nabla$ of the normal bundle
$\Gamma(T\Sigma^\perp)$ defined by 
\be
\tilde\nabla_X Z = \le\{ \matrix{
\displaystyle{ \pi [X,Z]\ ,\qquad X\in \Gamma(T\Sigma)}\cr
 \displaystyle{ \pi (\nabla_X Z)\ ,\qquad X\in \Gamma(T\Sigma^\perp)}
}\ri.\ ,
\ee
where $Z\in \Gamma(T\Sigma^\perp)$ is a normal section and $\pi:TM\to
T\Sigma^\perp$ is the natural projection. 
A classical result of \cite{kn:to} ensures that $\tilde\nabla$ is metric if
and only if ${\cal F}$ is Riemannian and $g$ is bundle--like.\\
Writing $Z=a^i\O_i$ and letting $A^i= \hbox{grad} (a^i)$, we have for $X\in
\Gamma(T\Sigma)$
\beaq
 [X,Z]&& = <X,A^i>\O_i + a^i [X,\O_i]= \cr
&& = <X,A^i>\O_i + a^i\le(\nabla_X\O_i-\nabla_{\O_i} X\ri)=\cr
&& =  <X,A^i>\O_i - a^i\le(\k\o_i X+\nabla_{\O_i} X\ri)\ .
\eeaq
On the other hand $<\nabla_{\O_i}X,\O_j> = -<X,\nabla_{\O_i}\O_j> = \k
\o_j<X,\O_i>=0$, namely $\nabla_{\O_i}X\in T\Sigma$.
Therefore we find 
\be
\tilde \nabla_XZ = \pi [X,Z] = <X,A^i>\O_i\ .
\ee
We not take $X\in \Gamma (T\Sigma^\perp)$ and compute
\be
\nabla_XZ = <X,A^i> \O_i + a^i\nabla_X\O_i= <X,A^i>\O_i -\k a^i\o_i X\ ,
\ee 
from which we obtain 
\be
\tilde\nabla_XZ = \pi(\nabla_XZ) = \nabla_XZ\ .
\ee
Summarizing the expression for the connection we have
\be
\tilde \nabla_X Z = \le\{
\matrix{
\displaystyle{ <X,A^i> \O_i, \qquad X\in T\Sigma}\cr
\displaystyle{ \nabla_X Z , \qquad X\in T\Sigma^\perp}
}\ri. 
\ee
The proof that $\tilde \nabla$ is a metric connection is a
straightforward computation of $X\|Z\|^2$ and $2<\tilde \nabla_XZ,Z>$ and
subsequent comparison.
$\rule{5pt}{8pt}$ \par\vskip 3pt
We conclude with a description of the infinitesimal automorphisms of $\cal F$
and ${\cal F}^\perp \equiv 
{\cal S}$ (the distribution spanned by $\O_1,...,\O_m$).\\
For $X\in \Gamma({\cal F})$ and $Z\in \Gamma({\cal S})$ we obtain
\be
[Z,X]=a^i\nabla_{\O_i} X-<X,A^i> \O_i + \k a^i\o_i X\ .
\ee
This shows that $Z$ is an infinitesimal automorphism of $\cal F$ (namely
$[Z,X]\in {\cal F}$) iff $<X,A^i>=0$. Thus
\bp
The infinitesimal automorphisms of $\cal F$ are the sections of $\Gamma ({\cal F})$
and combinations $a^i\O_i$ whose coefficients are constant on the leaves of
$\cal F$.
\ep
Finally, $X$ is an infinitesimal automorphism of ${\cal F}^\perp= {\cal S}$
iff $\nabla_{\O_i} X+ \k \o_i X=0$, namely iff $[X,\O_i]=0$.
\bp
The infinitesimal automorphisms of $\cal S$ are the sections of $\Gamma({\cal
S})$ and all sections of $\cal F$ which commute with every $\O_i$.
\ep

\end{document}